\documentclass[a4paper,11pt,reqno]{amsart}

\usepackage{supertabular}

\usepackage{array}

\usepackage[all]{xy}
\usepackage[latin1]{inputenc}

\def\theoremnumberstyle{changesection} 
\def\proofname{Proof}

\usepackage{michael}

\theoremstyle{break-italic}
\swapnumbers
\newtheorem{Theorem}{Theorem}

\numberwithin{equation}{section} 
\allowdisplaybreaks[3]   

\title{Quasi-Homogeneous Linear Systems on $¶^2$ with Base Points of Multiplicity $6$}
\author{Michael Kunte}
\address{Fachbereich Mathematik, Technische Universität
  Kaiserslautern, Erwin-Schrödinger-Straße, 67663 Kaiserslautern}
\email{kunte@mathematik.uni-kl.de}
\begin{document}

\maketitle

\begin{abstract}
 In this paper we prove the Harbourne-Hirschowitz conjecture for
 quasi-homogeneous linear systems of multiplicity $6$ on $¶^2$. For
 the proof we use the degeneration of the plane by Ciliberto
 and Miranda and results by Laface, Seibert, Ugaglia and Yang. As an
 application we derive a classification of the special systems of
 multiplicity $6$. 
\end{abstract}

\section{Introduction}

A classical problem in algebraic geometry is the dimensionality problem
for plane curves, which can be formulated as follows. Given finitely 
many general points of the projective plane with assigned multiplicities 
and a number $d$, determine the dimension of the linear system of curves 
of degree $d$ having at the given points at least the assigned multiplicities.
More precisely, the problem is to classify all systems which fail to
have the expected dimension (see \cite{CilIII} for some remarks on the 
history of this problem and its geometric meaning). Harbourne and Hirschowitz 
conjecture that these special systems are precisely the $(-1)$-special systems. In this paper, we give a complete list of the $(-1)$-special systems in the case in which the assigned
multiplicity is $6$ at all but one of the given points. Our main
result is the proof of the  Harbourne-Hirschowitz conjecture in this case.

We proceed along the following lines.
In Section 2 we introduce the necessary notation and give
a precise statement of the Harbourne-Hirschowitz conjecture.
In Section 3 we present a list of the $(-1)$-special linear systems
in our case. Its completeness is proved in Section 4.
In Section 5 we review the degeneration of the plane by
Ciliberto and Miranda. This method is the key tool in our proof 
of the main result which is given in the final two sections. 

\section{The Harbourne-Hirschowitz conjecture}

We work over the complex numbers and choose $n+1$ general points
$p_0,p_1,\ldots,p_n$ in $¶^2$, the projective plane over that
field.

\begin{notation}
We write $\kl = \kl(d,m_0,m_1,\ldots,m_n) \subset
¶(\Gamma(¶^2,\ko_{¶^2}(d)))$ for the linear system of all curves of
degree $d$ in $¶^2$ having multiplicity at least $m_i$ at $p_i$
for all $i$. 
We denote by $\ell(\kl)$ its projective dimension.

\end{notation}

Let $¶'$ be the blow-up of $¶^2$ at $p_0,p_1,\ldots,p_n$. By $H$
we denote the pull-back of a line in $¶^2$ and by $E_i$ the
exceptional divisor over $p_i$. The dimension of $\kl$ is the same
as the dimension of $|D|$ on $¶'$ with $D =
dH-m_0E_0-m_1E_1-\ldots-m_nE_n$. Using cohomology on $¶'$, we have
\begin{displaymath}
       \ell(\kl) = h^0(\ko_{¶'}(D)) - 1.
\end{displaymath}
Therefore we have by Riemann-Roch
\begin{displaymath}
      \ell(\kl) = \frac{D.(D-K_{¶'})}{2} + h^1(\ko_{¶'}(D)) - h^2(\ko_{¶'}(D)) + \chi(\ko_{¶'}) - 1
\end{displaymath}
($K_{¶'}$ denotes the canonical divisor on $¶'$). Since the
arithmetic genus of $¶'$ is zero, Serre duality implies
\begin{displaymath}
  \ell(\kl) = \frac{D.(D-K_{¶'})}{2} + h^1(\ko_{¶'}(D)).
\end{displaymath}

\begin{definition}

We define the \textit{virtual dimension} $v(\kl)$ of $\kl$ as follows:
   \begin{displaymath}
     v(\kl) = \frac{D.(D-K_{¶'})}{2}.
   \end{displaymath}
We define the \textit{expected dimension} to be
   \begin{displaymath}
     e(\kl) = \max \{-1,v(\kl) \}.
   \end{displaymath}
As $v(\kl) = \frac{d(d+3)}{2} - \sum_{i=0}^{n} \frac{m_i(m_i+1)}{2}$, one sees that the expected dimension is the one we obtain if all conditions imposed on the base points are independent.

We define $\kl$ to be \textit{special} or \textit{non-regular} if
  \begin{displaymath}
    \ell(\kl) > e(\kl),
  \end{displaymath}
otherwise we call $\kl$ \textit{non-special} or \textit{regular}.
\end{definition}

We recall some definitions from \cite{CilI}:

\begin{definition}[$(-1)$-special systems]
    Let $\ka$ in $¶^2$ be an irreducible curve such that its strict
    transform $\tilde \ka$ in $¶'$ is rational and smooth. Then $\ka$ is a \textit{(-1)-curve} if the self-intersection number
             \begin{displaymath}
                    \tilde \ka^{2} = -1.
             \end{displaymath}
By $\kl.\ka$ we denote the intersection number $D.\tilde \ka$ on $¶'$.

The linear system $\kl$ is called \textit{(-1)-special} if
  \begin{itemize}

     \item there exist $\ka_1,\ldots,\ka_t$ $(-1)$-curves with \mbox{$\kl.\ka_i =
  -n_i$} such that $n_i \geq 1$ for all $i$,

    \item there is an $j$ with $n_j \geq 2$ and

    \item the residual system \mbox{$\km = \kl - \sum_{i=0}^{t} n_i \ka_i$} has $v(\km) \geq 0$.

  \end{itemize}
\end{definition}

The main conjecture can be formulated as follows:

\begin{conjecture}[Harbourne-Hirschowitz]
     A linear system $\kl=\kl(d,m_0,m_1,\ldots,m_n)$ is special if and only if it is $(-1)$-special.
\end{conjecture}

It is easy to see that a $(-1)$-special system $\kl$ is
special because
\begin{displaymath}
v(\kl) = \frac{\kl.(\kl-K_{¶'})}{2} = \frac{(\km
  + n \ka).(\km + n \ka - K_{¶'})}{2}.
\end{displaymath}
Since $\ka.K_{¶'} = -1$ by the rationality of $\tilde \ka$, this implies
\begin{displaymath}
 v(\kl) = v(\km) + \frac{- n^2 + n}{2}
\leq \ell(\kl) + \frac{ - n^2+n}{2}.
\end{displaymath}
Therefore the opposite direction of the Harbourne-Hirschowitz
conjecture is the non-trivial one. It states that every
special system $\kl$ has fixed multiple $(-1)$-curves. Proving
the conjecture leads to an answer of the dimensionality problem.

\begin{remark}
We give a list of results on the conjecture. In fact we use all of
them in several ways for the proof of our main theorem.

We write $\kl = \kl(d,m_0^{b_0},m_1^{b_1},\ldots,m_r^{b_r})$ if $\kl$ has precisely $b_i$ base
 points of multiplicity $m_i$ for $i=0,\ldots,r$. With this notation the conjecture holds if
\begin{itemize}

  \item $b_0 + \ldots + b_r \leq 9$ \cite{Har},

  \item $\kl = \kl(d,m^n)$ (call it \textit{homogeneous of
      multiplicity m}) and $m \leq 12$ \cite{CilII},

  \item $\kl = \kl(d,m_0,m^n)$ (call it \textit{quasi-homogeneous
      of multiplicity m}) and $m \leq 3$ \cite{CilI},

  \item $\kl = \kl(d,m_0,4^n)$ \cite{Seib} and \cite{LafII},

  \item $\kl = \kl(d,m_0,5^n)$ \cite{Laf} or

  \item all multiplicities are bounded by $6$, i.e. $m_i \leq 6$ for $i=0,1,\ldots,n$ \cite{Yang}.

\end{itemize}

\end{remark}

\section{Main Results}

Our main result is a proof of the Harbourne-Hirschowitz conjecture in the quasi-homogeneous case of multiplicity $6$:

\begin{Theorem}[Main Theorem]
\label{main-theorem}
A system $\kl(d,m_0,6^n)$ is special if and only if it is $(-1)$-special.
\end{Theorem}

We give the proof within an extra section. For the proof we need the
following classification:

\begin{Theorem}[Classification of $(-1)$-special systems $\kl(d,m_0,6^n)$]

\label{classification}

The following is a complete list of all $(-1)$-special
systems $\kl(d,m_0,6^n)$.

\tablefirsthead{ $d-m_0$ & system & $v(\kl)$ & $\ell(\kl)$ & \\

   \hline
\\}

\tablehead{ $d-m_0$ & system & $v(\kl)$ & $\ell(\kl)$ & \\

   \hline
\\}

\begin{supertabular*}{\textwidth}%
{@{}c@{\extracolsep{\fill}}l@{\extracolsep{\fill}}l@{\extracolsep{\fill}}l@{\extracolsep{\fill}}l@{}}

$0$ &$ \mathcal{L}(d,d,6^{n})$ & $-21n+d$ & $-6n+d$ & $d \geq 6n \geq 6$
\\
$1$ &$ \mathcal{L}(d,d-1,6^{n})$ & $-21n+2d$ & $-11n+2d$ & $d \geq \frac{11}{2}n \geq \frac{11}{2}$
\\
$2$ &$ \mathcal{L}(10e,10e-2,6^{2e})$ & $-12e-1$ & $0$ & $e \geq 1$\\

 &$ \mathcal{L}(d,d-2,6^{n})$ & $-21n+3d-1$ & $-15n+3d-1$ & $d \geq \frac{1+15n}{3} \geq \frac{16}{3}$
\\
$3$&$ \mathcal{L}(9e,9e-3,6^{2e})$ & $-6e-3$ & $0$ & $e \geq 1$
\\
&$ \mathcal{L}(9e+1,9e-2,6^{2e})$ & $-6e+1$ & $2$ & $e \geq 1$
\\
&$ \mathcal{L}(d,d-3,6^{n})$ & $-21n+4d-3$ & $\geq -18n+4d-3$ & $d \geq \frac{18n+3}{4} \geq \frac{21}{4}$
\\
&&& \scriptsize $=$ if $d \neq \frac{9n}{2}+1$ or $n$ odd&
\\
$4$&$ \mathcal{L}(8e,8e-4,6^{2e})$ & $-2e-6$ & $0$ & $e \geq 1$
\\
&$ \mathcal{L}(8e+1,8e-3,6^{2e})$ & $-2e-1$ & $2$ & $e \geq 1$
\\
&$ \mathcal{L}(8e+2,8e-2,6^{2e})$ & $-2e+4$ & $5$ & $e \geq 1$
\\
&$ \mathcal{L}(d,d-4,6^{n})$ & $-21n+5d-6$ & $\geq -20n+5d-6$ & $d \geq \frac{20n+6}{5} \geq \frac{26}{5}$
\\
&&& \scriptsize $=$ if $d \neq 4n+2$ or $n$ odd&
\\
$5$&$ \mathcal{L}(7e,7e-5,6^{2e})$ & $-10$ & $0$ & $e \geq 1$
\\
&$ \mathcal{L}(7e+1,7e-4,6^{2e})$ & $-4$ & $2$ & $e \geq 1$
\\
&$ \mathcal{L}(7e+2,7e-3,6^{2e})$ & $2$ & $5$ & $e \geq 1$
\\
&$ \mathcal{L}(7e+3,7e-2,6^{2e})$ & $8$ & $9$ & $e \geq 1$
\\
$6$&$ \mathcal{L}(6e,6e-6,6^{2e})$ & $-15$ & $0$ & $e \geq 1$
\\
&$ \mathcal{L}(6e+1,6e-5,6^{2e})$ & $-8$ & $2$ & $e \geq 1$
\\
&$ \mathcal{L}(6e+2,6e-4,6^{2e})$ & $-1$ & $5$ & $e \geq 1$
\\
&$ \mathcal{L}(6e+3,6e-3,6^{2e})$ & $6$ & $9$ & $e \geq 1$
\\
&$ \mathcal{L}(6e+4,6e-2,6^{2e})$ & $13$ & $14$ & $e \geq 1$           
\\
$7$&$ \mathcal{L}(5e+2,5e-5,6^{2e})$ & $-2e-5$ & $-2e+5$ & $2 \geq e \geq 1$
\\
&$ \mathcal{L}(5e+3,5e-4,6^{2e})$ & $-2e+3$ & $-2e+9$ & $4 \geq e \geq 1$
\\
&$ \mathcal{L}(5e+4,5e-3,6^{2e})$ & $-2e+11$ & $-2e+14$ & $7 \geq e \geq 1$
\\
&$ \mathcal{L}(5e+5,5e-2,6^{2e})$ & $-2e+19$ & $-2e+20$ & $10 \geq e \geq 1$
\\
$8$&$ \mathcal{L}(4e+4,4e-4,6^{2e})$ & $-6e+8$ & $-6e+14$ & $2 \geq e \geq 1$
\\
&$ \mathcal{L}(4e+5,4e-3,6^{2e})$ & $-6e+17$ & $-6e+20$ & $2 \geq e \geq 1$
\\
&$ \mathcal{L}(4e+6,4e-2,6^{2e})$ & $-6e+26$ & $-6e+27$ & $4 \geq e
\geq 1$
\\
&$ \mathcal{L}(10,2,6^{3})$ & $-1$ & $2$ &
\\
&$ \mathcal{L}(24,16,6^{9})$ & $-1$ & $0$ &
\\
$9$&$ \mathcal{L}(3e+6,3e-3,6^{2e})$ & $-12e+24$ & $-12e+27$ & $2 \geq e \geq 1$
\\
&$ \mathcal{L}(3e+7,3e-2,6^{2e})$ & $-12e+34$ & $-12e+35$ & $2 \geq e \geq 1$
\\
&$ \mathcal{L}(9,0,6^{3})$ & $-9$ & $0$ &
\\
&$ \mathcal{L}(10,1,6^{3})$ & $1$ & $4$ &
\\
&$ \mathcal{L}(14,5,6^{5})$ & $-1$ & $0$ &
\\
&$ \mathcal{L}(18,9,6^{7})$ & $-3$ & $0$ &
\\
$10$&$ \mathcal{L}(2e+8,2e-2,6^{2e})$ & $-20e+43$ & $-20e+44$ & $2 \geq e \geq 1$
\\
&$ \mathcal{L}(10,0,6^{3})$ & $2$ & $5$ &
\\
&$ \mathcal{L}(14,4,6^{5})$ & $4$ & $5$ &
\\
$11$&$ \mathcal{L}(13,2,6^{5})$ & $-4$ & $2$ &
\\
&$ \mathcal{L}(14,3,6^{5})$ & $8$ & $9$ &
\\
$12$&$\mathcal{L}(12,0,6^{5})$ & $-15$ & $0$ &
\\
&$ \mathcal{L}(13,1,6^{5})$ & $-2$ & $4$ &
\\
&$ \mathcal{L}(14,2,6^{5})$ & $11$ & $12$ &
\\
$13$&$\mathcal{L}(13,0,6^{5})$ & $-1$ & $5$ &
\\
&$ \mathcal{L}(14,1,6^{5})$ & $13$ & $14$ &
\\
$14$&$\mathcal{L}(14,0,6^{5})$ & $14$ & $15$&
\\
\end{supertabular*}

\end{Theorem}

\section{The Classification}

In the paper \cite{CilI} of Ciliberto and Miranda a lot of
classification work has been done which we can apply to our
problem. Ciliberto and Miranda introduced two notions which we recall now to use their results.

Let $\kl$ be a linear system of plane curves with general multiple
base points as above. Then 
$\kl$ is a \textit{quasi-homogeneous (-1)-class} if $\kl = \kl(d,m_0,m^n)$, on $¶'$ the
self-intersection number $\kl.\kl = -1$ and the arithmetic genus
      \begin{displaymath}
           g_{\kl} = \frac{\kl^2+\kl.K_{¶'}}{2} +1 = 0.
      \end{displaymath}
 As $v(\kl) = \kl^2 - g_{\kl} +1$, these systems are never
 empty.

In this case, if $\ka$ is a $(-1)$-curve such that $\ka \in
\kl$ then by $\kl.\ka = -1$ and the irreducibility of $\ka$, we
have $\kl = \{\ka\}$. So we can identify $(-1)$-curves and quasi-homogeneous $(-1)$-classes and write
$\ka = \kl$. Ciliberto and Miranda proved that such a $(-1)$-curve
exists up to $m \leq 6$. Hence a numerical classification of these systems
gives a classification for all quasi-homogeneous
$(-1)$-curves up to multiplicity $m = 6$. Such a classification is
given in \cite{CilI}.

Now we consider the following phenomenon: Let $\kl =
\kl(d,m_0,m^n)$ be a quasi-homogeneous linear system and
$\ka$ a $(-1)$-curve such that $\ka =
\kl(\delta,\mu_0,\mu_1,\ldots,\mu_n)$ and $\kl.\ka \leq -2$. Let
$\Perm_n$ be the permutation group on $n$ letters and let $\sigma
\in \Perm_n$. We define $\ka_{\sigma} =
\kl(\delta,\mu_0,\mu_{\sigma(1)},\ldots,\mu_{\sigma(n)})$. Then,
as $\ka$ is a $(-1)$-curve, it follows that $\ka_{\sigma}$
is again a $(-1)$-curve. As $\kl$ is quasi-homogeneous we have again $\kl.\ka_{\sigma} \leq
-2$. Therefore we can construct a composition of
$(-1)$-curves, which split off the system $\kl$. We define
the set $A \subset \Perm_n$ to be maximal such that all
$\ka_{\sigma}$ with $\sigma \in A$ are pairwise different. Then we
define a new plane curve $\ka_{tot}= \sum_{\sigma \in A}
\ka_{\sigma}$ (see \cite{Laf}).

We call a linear system $\kl'=\kl(d,m_o,m_1,\ldots,m_n)$ as above a \textit{quasi-homogeneous (-1)-configuration} if $\ka_{\tot}$ is a generic element in $\kl'$. We note that $\kl'$ is by construction quasi-homogeneous (if $k = |A|$ then there exists a $\mu'$ such that $\kl' = \kl(k\delta,k\mu_0,\mu'^n)$).

\begin{lemma}[splitting-off Lemma]

\label{splitting-off}

Let $\kl = \kl(d,m_0,m^n)$. Then every $(-1)$-curve $\ka$
with $\kl.\ka \leq -2$ is of one of the following types (We have
listed the associated quasi-homogeneous compound
$(-1)$-configurations, too.):

   \begin{tabular}{ll}
          $ \ka = \kl(\delta,\mu_0,\mu_1^n)$ & \\
          $ \ka = \kl(\delta,\mu_0,\mu_2 - 1,\mu_2^{n-1})$& $\ka_{tot} = \kl(n\delta,n \mu_0, (n\mu_2 - 1)^{n})$\\
          $ \ka = \kl(\delta,\mu_0,\mu_2 + 1,\mu_2^{n-1})$& $\ka_{tot} = \kl(n\delta,n \mu_0, (n\mu_2 + 1)^{n})$\\
   \end{tabular}
\end{lemma}

\begin{proof}
First one proves that strict transforms of different $\ka_{\sigma}
\neq \ka_{\sigma'}$ cannot meet positively on $¶'$. This is the
case as otherwise one sees, by the Riemann-Roch theorem on $¶'$,
that the sum of these moves in a linear system of positive
dimension, which is a contradiction to being a fixed part of
$\kl$. This implies that all the different $\ka_{\sigma}$ are
linearly independent in $\Pic(¶')$. Let the $\mu_1,\ldots,\mu_n$ occur in sets
of size $k_1 \leq \ldots \leq k_s$. As $\rank \Pic(¶') = n+2$ we see
by combinatorial reasons that for the $\frac{n!}{k_1!\cdots k_s!}$
different $(-1)$-curves $\ka_{\sigma}$ only the
possibilities
\begin{displaymath}
\begin{array}{lc}
   s = 1,\; k_1 = n & \mbox{or} \\
   s=2 , \;k_1 = 1,\; k_2 = n-1&
\end{array}
\end{displaymath}
can occur. That means we have at most three different
multiplicities $\mu_0,\mu_1$ and $\mu_2$. 

Moreover we have the equations $\ka.\ka = -1$ and $\ka.\ka_{\sigma} = 0$ on $¶'$. That gives $\ka.\ka-\ka.\ka_{\sigma} = -1$ which is equivalent to $(\mu_1-\mu_2)^2 = 1$ (see \cite{CilI}).
\end{proof}

For the purpose of classifying the systems $\kl(d,m_0,6^n)$ we
need a complete list of all $(-1)$-curves which might split
off such systems two times. These $(-1)$-curves can not have
higher multiplicities than $3$ at the points $p_1,\ldots,p_n$. We
obtain the following result:

\begin{lemma}[classification of $(-1)$-curves]

\label{class-einscurves}

All $(-1)$-curves $\ka$ and quasi-homogeneous
$(-1)$-configurations $\ka_{\tot}$ up to multiplicity $3$ in the
points $p_1,\ldots,p_n$ which might split off a quasi-homogeneous
system $\kl = \kl(d,m_0,6^n)$ are elements of the systems in the
following list (see \cite{Laf}):

\begin{tabular}{ll}
not compound & compound
\\
$\mathcal{L}(2,0,1^{5})$&
\\
$\mathcal{L}(e,e-1,1^{2e})\quad e \geq 1$&
\\
$\mathcal{L}(1,1,1^{1})$&$\mathcal{L}(n,n,1^{n}) \quad n \geq 2$
\\
$\kl(1,0,1^2)$&$\mathcal{L}(3,0,2^{3})$
\\
$\mathcal{L}(6,3,2^{7})$&
\\
$\mathcal{L}(12,8,3^{9})$&

\end{tabular}

In particular, all the $(-1)$-curves are quasi-homogeneous.

\end{lemma}

\begin{proof}
We refer to \cite[Example 5.1]{CilI} for the proof of a list of all quasi-homogeneous $(-1)$-classes up to multiplicity $3$. 
In \cite[Example 5.15]{CilI} is given a complete list of all quasi-homogeneous $(-1)$-configurations up to multiplicity $3$. Using this two lists and Lemma \ref{splitting-off} gives this result.
\end{proof}

Now we give the proof of the classification theorem of all
$(-1)$-special systems of the form $\kl(d,m_0,6^n)$.

\begin{proof}[Proof of Theorem \ref{classification}]

In lemma \ref{class-einscurves} we have seen the possible cases for $(-1)$-curves which might split off $\kl(d,m_0,6^n)$. Now we have to consider all these cases.
To be a little bit faster we proceed along the following algorithm (see \cite{Laf}):

We go through all possible combinations of these $(-1)$-curves step by step.

First step:
If we find a $(-1)$-curve or a $(-1)$-configuration $\ka$ such that
\begin{displaymath}
  \kl.\ka = -\mu \leq -2,
\end{displaymath}
then we split off the fixed part and define $\km = \kl -
\mu\cdot\ka$.

Second step:
Let $\km'$ be the residual system of $\km$ obtained by splitting off all possible $(-1)$-curves. By the definition of $(-1)$-special systems we have to verify that $v(\km') \geq 0$. We notice that the systems $\km$ are quasi-homogeneous of multiplicity $\leq 4$ by lemma \ref{class-einscurves}. Therefore we can use the results of \cite{CilI} and \cite{Seib}.

\renewcommand{\labelitemii}{}

\begin{itemize}

\item  \framebox{\textbf{$\mathbf{\kl = \km + \mu \cdot \ka}$,
      $\mathbf{v(\km)\geq 0}$ and $\mathbf{\km.\ka = 0}$}}

\begin{enumerate}

\medskip

\item

\textbf{$\mathbf{\ka = \kl(2,0,1^5)}$ and $\mathbf{\kl = \kl(d,m_0,6^5)}$}

\medskip

This gives  $\km = \kl(d-2n,m_0,(6-\mu)^5)$ and $\km.\ka = 0$ gives $d = \frac{30-\mu}{2}$.
\medskip

If \ul{$\mu = 2$}  $\Lra$  $d=14$ and we get

\begin{itemize}

\item $m_0 = 0$ and $v(\km) = 15$ with $\km = \kl(10,0,4^5)$, non-special by \cite{Seib}

\item $m_0 = 1$ and $v(\km) = 14$ with $\km = \kl(10,1,4^5)$\mbox{,       } $''$

\item $m_0 = 2$ and $v(\km) = 12$ with $\km = \kl(10,2,4^5)$\mbox{,       } $''$

\item $m_0 = 3$ and $v(\km) = 9$ with $\km = \kl(10,3,4^5)$\mbox{,       }  $''$

\item $m_0 = 4$ and $v(\km) = 5$ with $\km = \kl(10,4,4^5)$\mbox{,       }  $''$

\item $m_0 = 5$ and $v(\km) = 0$ with $\km =  \kl(10,5,4^5)$\mbox{,       } $''$

\end{itemize}

\medskip

\ul{$\mu = 3$} is not possible because of $\km.\ka = 0$.

\medskip
If \ul{$\mu = 4$}  $\Lra$  $d=13$ and we conclude

\begin{itemize}

\item $m_0 = 0$ and $v(\km) = 5$ with $\km = \kl(7,0,2^5)$, non-special by \cite{CilI}
\item $m_0 = 1$ and $v(\km) = 4$ with $\km = \kl(7,1,2^5)$ \mbox{,       } $''$
\item $m_0 = 2$ and $v(\km) = 2$ with $\km = \kl(7,2,2^5)$ \mbox{,       } $''$
\item $m_0 = 3$ and $v(\km) = -1$
\end{itemize}

\medskip
\ul{$\mu=5$} is not possible because of $\km.\ka = 0$.

\medskip

From \ul{$\mu = 6$} $\Lra$  $d=12$ and $m_0=0$, $v(\km)=0$ for $\km=\kl(0,0)$.

\medskip

\item

\textbf{$\mathbf{\ka = \kl(e,e-1,1^{2e})}$ $\mathbf{e\geq 1}$ and $\mathbf{ \kl = \kl(d,m_0,6^{2e})}$}

\medskip

Then follows $\km = \kl(d-\mu\cdot e,m_0-\mu\cdot e+\mu,(6-\mu)^{2e})$ and $\km.\ka = 0$ gives \mbox{$-e\cdot m_0+e\cdot d -12e+m_0+\mu =0$}
$\Lra$  $m_0 > d-12$. $v(\km) \geq 0$ gives $d \geq m_0+\mu-2$.

\medskip
\ul{$\mu = 2$}  $\Lra$  $d \geq m_0 > d-12$

\begin{tabular}{ccc}

\hline
$m_0$ & $v(\km) \leq -1$ and $\km$ non-special&\\
\hline

$d$&   $''$&\\
$d-1$ & $''$&\\

\end{tabular}

\begin{tabular}{ccp{6cm}}

\hline
$m_0$& from $\km.\ka = 0 \Ra d $ & residual system\\
\hline

$d-2$&$10e$ & $\km= \kl(8e,8e,4^{2e})$ irregular by \Seib\space $\Ra$ non-empty \\

$d-3$&$9e+1$ & $\km = \kl(7e+1,7e,4^{2e})$ irregular by \Seib\space $\Ra$ non-empty\\
$d-4$&$8e+2$& $\km=\kl(6e+2,6e,4^{2e})$ irregular by \Seib\space $\Ra$ non-empty\\
$d-5$&$7e+3$& $\km=\kl(5e+3,5e,4^{2e})$ regular by \Seib \space and \mbox{$v(\km) = 9$}\\
$d-6$&$6e+4$& $\km=\kl(4e+4,4e,4^{2e})$ regular by \Seib \space and \mbox{$v(\km) = 14$}\\
$d-7$&$5e+5$& $\km=\kl(3e+5,3e,4^{2e})$ regular by \Seib \space and \mbox{$v(\km) = -2e+20$}\\
$d-8$&$4e+6$& $\km=\kl(2e+6,2e,4^{2e})$ regular by \Seib \space and \mbox{$v(\km) = -6e+27$}\\
$d-9$&$3e+7$& $\km=\kl(e+7,e,4^{2e})$ regular by \Seib \space and \mbox{$v(\km) = -12e+35$}\\
$d-10$&$2e+8$& $\km=\kl(8,0,4^{2e})$ regular by \Seib \space and \mbox{$v(\km) = -20e+44$}\\
$d-11$&$e+9$& $\Ra$ $m_0 \leq -1$ not possible\\
\end{tabular}

\medskip

\ul{$\mu = 3$}  $\Lra$  $d-1 \geq m_0 > d-12$

\begin{tabular}{ccc}

\hline
$m_0$ & $v(\km) \leq -1$ and $\km$ non-special&\\
\hline

$d-1$&$''$&\\
$d-2$ & $''$&\\
\end{tabular}

\tablefirsthead{\hline
                $m_0$ & from $\km.\ka = 0 \Ra d$ & residual system\\
                \hline}
\tablehead{\hline 
           $m_0$ & from $\km.\ka = 0 \Ra d$ & residual system\\
           \hline}

\begin{supertabular}{ccp{6cm}}
$d-3$ & $9e$& $\km=\kl(6e,6e,3^{2e})$ irregular by \CilI\space and \mbox{$e(\km)=0$}\\
$d-4$&$8e+1$& $\km=\kl(5e+1,5e,3^{2e})$ irregular by \CilI \space and \mbox{$e(\km) = 2$}\\
$d-5$&$7e+2$& $\km=\kl(4e+2,4e,3^{2e})$ regular by \CilI \space and \mbox{$v(\km) = 5$}\\
$d-6$&$6e+3$& $\km=\kl(3e+3,3e,3^{2e})$ regular by \CilI \space and \mbox{$v(\km) = 9$}\\
$d-7$&$5e+4$& $\km=\kl(2e+4,2e,3^{2e})$ regular by \CilI \space and \mbox{$v(\km) = -2e+14$}\\
$d-8$&$4e+5$& $\km=\kl(e+5,e,3^{2e})$ regular by \CilI \space and \mbox{$v(\km) = -6e+20$}\\
$d-9$&$3e+6$& $\km=\kl(6,0,3^{2e})$ regular by \CilI  \space and \mbox{$v(\km) = -12e+27$}\\
$d-10$&$2e+7$& $\Ra$ $m_0 \leq -1$ not possible\\
\end{supertabular}

\medskip

\ul{$\mu = 4$}  $\Lra$  $d-2 \geq m_0 > d-12$

\begin{tabular}{ccp{6cm}}

\hline
$m_0$ & from $\km.\ka = 0 \Ra d$ & residual system \\
\hline

$d-2$&$10e-2$& $\km = \kl(6e-2,6e,2^{2e})$ empty \\
$d-3$&$9e-1$& $\km=\kl(5e-1,5e,2^{2e})$ empty\\
$d-4$&$8e$& $\km=\kl(4e,4e,2^{2e})$ irregular by \CilI \space and \mbox{$e(\km) = 0$}\\
$d-5$&$7e+1$& $\km=\kl(3e+1,3e,2^{2e})$ regular by \CilI \space and \mbox{$v(\km) = 2$}\\
$d-6$&$6e+2$& $\km=\kl(2e+2,2e,2^{2e})$ regular by \CilI \space and \mbox{$v(\km) = 5$}\\
$d-7$&$5e+3$& $\km=\kl(e+3,e,2^{2e})$ regular by \CilI  \space and \mbox{$v(\km)=-2e+9$}\\
$d-8$&$4e+4$& $\km=\kl(4,0,2^{2e})$ regular by \CilI  \space and \mbox{$v(\km) = -6e+14$}\\
$d-9$&$3e+5$& $\Ra$ $m_0 \leq -1$ not possible\\

\end{tabular}

\medskip

\ul{$\mu = 5$}  $\Lra$  $d-3 \geq m_0 > d-12$

\begin{tabular}{ccp{6cm}}
\hline
$m_0$ & from $\km.\ka=0 \Ra d$ &residual system \\
\hline

$d-3$&$9e-2$& $\km=\kl(4e-2,4e,1^{2e})$ empty\\
$d-4$&$8e-1$& $\km=\kl(3e-1,3e,1^{2e})$ empty\\
$d-5$&$7e$& $\km=\kl(2e,2e,1^{2e})$ regular by \CilI \space and \mbox{$v(\km) = 0$}\\
$d-6$&$6e+1$& $\km=\kl(e+1,e,1^{2e})$ regular by \CilI \space and \mbox{$v(\km) = 2$}\\
$d-7$&$5e+2$& $\km=\kl(2,0,1^{2e})$ regular by \CilI  \space and \mbox{$v(\km)=-2e+5$}\\
$d-8$&$4e+3$& $\Ra$ $m_0 \leq -1$ not possible\\

\end{tabular}

\medskip

For \ul{$\mu = 6$} we have that $d-4 \geq m_0 > d-12$. Let $m_0=d-x$. From $\km.\ka=0$ $\Ra$ $d=(12-x)e+(x-6)$. We notice that $\km=\kl((6-x)e+(x-6),(6-x)e,0)$, which is regular. Taking into account that $v(\km) \leq -1$ for all $x \leq 5$ and $m_0 \leq -1$ for all $ x \geq 7$ we get the only case:

$m_0=d-6$ and $\km.\ka=0$ $\Ra$ $d=6e$ and $\km=\kl(0,0)$ is regular with $v(\km) = 0$.

\medskip

\item
\textbf{$\mathbf{\ka = \kl(e,e,1^e)}$ and $\mathbf{\kl = \kl(d,m_0,6^e)}$}

\medskip

This leads to $\km = \kl(d-\mu e,m_0-\mu e,(6-\mu)^e)$. $\km.\ka = 0$ gives $m_0 = d+\mu-6$.

If $\ul{\mu=2}$ then we get $m_0=d-4$, $\kl = \kl(d,d-4,6^e)$ and $\km=\kl(d-2e,d-4-2e,4^e)$. From $v(\km)=-20e+5d-6$ $\Lra$ $v(\km)\geq 0$ if $d\geq \frac{6+20e}{5}$. Furthermore $\km$ is irregular by \cite{Seib} and of higher dimension if

\begin{enumerate}

\item     $e=2f$ and $d=8f$
\item     $e=2f$ and $d=8f+1$
\item     $e=2f$ and $d=8f+2$.

\end{enumerate}

\medskip

If $\ul{\mu=3}$ then we get $m_0=d-3$, $\kl = \kl(d,d-3,6^e)$ and $\km=\kl(d-3e,d-3-3e,3^e)$. From $v(\km)=-18e+4d-3$ $\Lra$ $v(\km)\geq 0$ if $d\geq \frac{3+18e}{4}$. Further $\km$ is irregular by \CilI\space and of higher dimension if

\begin{enumerate}

\item $e=2f$, $d=9f$ and $e(\km)=0$ or
\item $e=2f$, $d=9f+1$ and $e(\km)=2$.

\end{enumerate} 

\medskip

If $\ul{\mu=4}$ then $m_0=d-2$, $\kl = \kl(d,d-2,6^e)$ and $\km=\kl(d-4e,d-2-4e,2^e)$. From $v(\km)=-15e+3d-1$ $\Lra$ $v(\km)\geq 0$ if $d\geq \frac{1+15e}{3}$. Further $\km$ is irregular by \CilI\space and of higher dimension if $e=2f$, $d=10f$ and $e(\km)=0$.

\medskip

If $\ul{\mu=5}$ then $m_0=d-1$, $\kl = \kl(d,d-1,6^e)$ and $\km=\kl(d-5e,d-1-5e,1^e)$. From $v(\km)=-11e+2d$ $\Lra$ $v(\km)\geq 0$ if $d\geq \frac{11e}{2}$. $\km$ is always regular by \CilI. 

\medskip

If $\ul{\mu=6}$ then $m_0=d$, $\kl = \kl(d,d,6^e)$ and $\km=\kl(d-6e,d-6e)$. $v(\km)=-6e+d$ $\Lra$ $v(\km)\geq 0$ if $d\geq 6e$. $\km$ is always regular. \\

\medskip

The following two cases are easier to compute because we have no further parameters in the \textit{(-1)-curves}.

\medskip
\item

\textbf{$\mathbf{\ka = \kl(6,3,2^7)}$ and $\mathbf{\kl = \kl(d,m_0,6^7)}$}

\medskip
This leads to $\km = \kl(d-6\mu,m_0-3\mu,(6-2\mu)^7)$. $\km.\ka = 0$ gives $m_0 = \frac{6d+\mu-84}{3}$. Therefore $\ul{\mu=3}$ is the only possible case: $\km=\kl(d-18,2d-36)$. To get $v(\km)\geq 0$ we need $d=18$. $\Lra$ $\kl=\kl(18,9,6^7)$.

\medskip

\item

\textbf{$\mathbf{\ka = \kl(3,0,2^3)}$ and $\mathbf{ \kl = \kl(d,m_0,6^3)}$}

\medskip

This leads to $\km = \kl(d-3\mu,m_0,(6-2\mu)^3)$. $\km.\ka = 0$ gives $d = 12-\mu$.

  \begin{tabular}{lp{12cm}}
     $\ul{\mu=2}$& We get $\kl=\kl(10,m_0,2^3)$ and $\km=\kl(4,m_0,2^3)$. From $v(\km) \geq 0$ we get $m_0 \in \{0,1,2\}$. All $\km$ are regular by \CilI.  \\
     $\ul{\mu=3}$& We get $\kl=\kl(9,m_0)$ and $\km=\kl(0,m_0)$. $\Lra$ $m_0=0$ and $v(\km) = 0$.\\
  \end{tabular}

\medskip

\item

  \textbf{$\mathbf{\ka = \kl(12,8,3^9)}$, $\mathbf{ \kl = \kl(d,m_0,6^9)}$ and $\mathbf{\mu=2}$} 

\medskip

This lead to $\Lra$ $\km=\kl(d-24,m_0-16)$, which is regular. From $\km.\ka = 0$ we get $m_0=\frac{3d-40}{2}$. Therefore $v(\km)\geq 0$ gives $d \in \{24,25\}$, but only $d=24$ and $m_0 =16$ is possible.

\end{enumerate}

\end{itemize}

\medskip

\item

\framebox{\textbf{$\mathbf{\kl = \km + 2 \cdot \ka_1 + 2 \cdot \ka_2}$,  $\mathbf{ v(\km)\geq 0}$,  $\mathbf{\km}$ non-special and $\mathbf{\km.\ka = 0}$}}

\begin{enumerate}

\item

\textbf{$\mathbf{\ka = \kl(\delta,\mu_0,1^n)}$ and $\mathbf{\ka_1.\ka_2 = 0}$}

\medskip

This leads to $\ka_1 = \kl(e,e-1,1^{2e})$ and $\ka_2 = \kl(2e,2e,1^{2e})$. Further we have $\kl=\kl(d,m_0,6^{2e})$ and $\km = \kl(d-6e,m_0-6e+2,2^{2e})$. From $\km.\ka_1 = 0$ and $\km.\ka_2 = 0$ we get $m_0 = d-4$ and $d=8e+2$. Therefore we have $\km=\kl(2e+2,2e,2^{2e})$, which is regular by \cite{CilI} and $v(\km)=5$.

\medskip

\item
\textbf{$\mathbf{\ka_1 = \kl(\delta_1,\mu_{0_1},1^n)}$ and $\mathbf{\ka_2 = \kl(\delta_2,\mu_{0_2},2^n)}$}

$\ka_1.\ka_2 = 0$ gives only the following possibilities:

\medskip

\begin{enumerate}
 \item $\ka_1 \in \kl(2,1,1^4)$ and $\ka_2 \in \kl(3,0,2^3)$,
 \item  $\ka_1 \in \kl(2,0,1^5)$ and $\ka_2 \in \kl(3,0,2^3)$ or
 \item  $\ka_1 \in \kl(2,2,1^2)$ and $\ka_2 \in \kl(3,0,2^3)$.
\end{enumerate}

(1) and (2) are not possible cases as these curves are elements of quasi-homogeneous systems based on a different number of points with equal multiplicities. It is not possible to find a suitable system $\kl(d,m_0,6^n)$. So let us focus on (3), where we see that it is equivalent to assume $\ka_2 \in \kl(3,2,2^2)$. In this case we see that $\kl=\kl(d,m_0,6^2)$ and $\km = \kl(d-10,m_0-8)$. From $\km.\ka_1 = 0$ and $\km.\ka_2 = 0$ we conclude that $d=10$ and $m_0=8$, that means we get the system $\kl=\kl(10,8,6^2)$.

\end{enumerate}

\medskip

\item

\framebox{\textbf{$ \mathbf{\kl = \km + 2 \cdot \ka_1 + 3 \cdot \ka_2}$,  $ \mathbf{v(\km)\geq 0}$ and $\mathbf{\km.\ka = 0}$}}

\framebox{\textbf{$\mathbf{\ka = \kl(\delta,\mu_0,1^n)}$ and $\mathbf{\ka_1.\ka_2 = 0}$}}

\begin{enumerate}

\item

\textbf{$\mathbf{\ka_1 = \kl(e,e-1,1^{2e})}$ $\&$ $\mathbf{\ka_2 = \kl(2e,2e,1^{2e})}$}

\medskip

Moreover we have $\kl=\kl(d,m_0,6^{2e})$ and $\km = \kl(d-8e,m_0-8e+2,1^{2e})$. From $\km.\ka_1 = 0$ and $\km.\ka_2 = 0$ we get $m_0 = d-3$ and $d=9e+1$. Therefore we have $\km=\kl(e+1,e,1^{2e})$ which is regular by \CilI\space and $v(\km)=2$.

\medskip

\item

\textbf{$\mathbf{\ka_1 = \kl(2e,2e,1^{2e})}$ $\&$ $\mathbf{\ka_2 =
    \kl(e,e-1,1^{2e})}$}

\medskip

Furthermore we have $\kl=\kl(d,m_0,6^{2e})$ and $\km = \kl(d-7e,m_0-7e+3,1^{2e})$. From $\km.\ka_1 = 0$ and $\km.\ka_2 = 0$ we get $m_0 = d-4$ and $d=8e+1$. Therefore we have $\km=\kl(e+1,e,1^{2e})$ which is regular by \CilI\space and $v(\km)=2$.

\end{enumerate}

\medskip

\item

 \framebox{\textbf{$\mathbf{\kl = \km + 2 \cdot \ka_1 + 4 \cdot \ka_2}$,  $\mathbf{ v(\km)\geq 0}$ and $\mathbf{\km.\ka = 0}$}}
\\
\framebox{\textbf{$\mathbf{\ka = \kl(\delta,\mu_0,1^n)}$ and $\mathbf{\ka_1.\ka_2 = 0}$}}

\bigskip

\begin{enumerate}

\item

\textbf{$\mathbf{\ka_1 = \kl(e,e-1,1^{2e})}$ and $\mathbf{\ka_2 = \kl(2e,2e,1^{2e})}$}

\medskip

Moreover we have $\kl=\kl(d,m_0,6^{2e})$ and $\km = \kl(d-10e,m_0-10e+2)$. From $\km.\ka_1 = 0$ and $\km.\ka_2 = 0$ we get $m_0 = d-2$ and $d=10e$. Therefore we get $\km=\kl(0,0)$ and $v(\km)=0$.

\medskip

\item

\textbf{$\mathbf{\ka_1 = \kl(2e,2e,1^{2e})}$ and $\mathbf{\ka_2 = \kl(e,e-1,1^{2e})}$}

\medskip

Furthermore we have $\kl=\kl(d,m_0,6^{2e})$ and $\km = \kl(d-8e,m_0-8e+4)$. From $\km.\ka_1 = 0$ and $\km.\ka_2 = 0$ we get $m_0 = d-4$ and $d=8e$. Therefore we have $\km=\kl(0,0)$ and $v(\km)=0$.

\end{enumerate}

\medskip

\item

 \framebox{\textbf{$\mathbf{\kl = \km + 2 \cdot \ka_1 + 2 \cdot \ka_2+ 2 \cdot \ka_3}$,  $\mathbf{ v(\km)\geq 0}$ and $\mathbf{\km.\ka = 0}$}}

\medskip

As $\ka_1 = \kl(e,e-1,1^{2e})$ and $\ka_2 = \kl(e,e,1^{e})$ are the only compound $(-1)$-configurations with multiplicity $m=1$ in $p_1,\ldots,p_n$ which have intersection multiplicity $= 0$. Therefore we are immediately in case 4.

\bigskip

This finally completes our proof of the classification theorem.
\end{proof}

\section{The Degeneration Method}

In this section we give a rough overview of the degeneration
 of the plane as introduced by Ciliberto and Miranda in
\cite{CilI}. We refer to this paper for further details. As in
every degeneration method the aim is to specialize the base points
of a system $\kl(d,m_0,m^n)$ in such a way that on the one hand
the dimension is easier to compute but on the other hand it does
not change.

At first we consider the geometric situation. Let $\Delta$ be a complex disc around the origin. We define $V = ¶^2 ×
\Delta$. Let $p_1: V \lra ¶^2$ and $p_2: V \lra \Delta$ be the projections. Now
we blow up a line $L$ in $V_0 = p_2^{-1}(0)$ ($f: X \lra V$) and obtain the following situation with $\pi_i = f \circ p_i$:

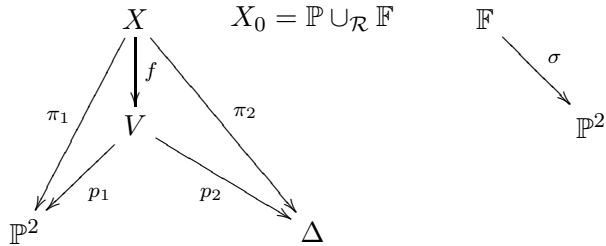
\xymatrix{&&&&&& X \ar[ddl]_{\pi_{1}} \ar[ddr]^{\pi_2} \ar[d]^{f} & X_0 = ¶\cup_{\kr} \F & \F \ar[dr]^{\sigma} &\\
          &&&&&& V \ar[dl]^{p_1} \ar[dr]_{p_2} & && ¶^2\\
          &&&&& ¶^2 & & \Delta &&}

Now $X_t = \pi_2^{-1}(t) \cong
¶^2$ for all $t \neq 0$. $X_0 = \pi_2^{-1}(0)$ is a union of two
surfaces, the strict transform of $V_0 \cong ¶^2$ (called $¶$) and the
exceptional divisor $\F = f^{-1}(L)$. $\F$ is isomorphic to the
blow-up of $¶^2$ in one point $p$ (here via $\sigma$). The surfaces are
glued together along the line $\kr$, which can be identified with $L$
in $¶$ and with the exceptional divisor $E = \sigma^{-1}(p)$ in $\F$.

As in \cite{CilI} we define $\ko_{X}(d) = \pi_1^{*}\ko_{¶^2}(d)$ and $\ko_{X}(d,k) = \ko_{X}(d) \otimes_{\ko_{X}} \ko_{X}(k¶)$. We set $\chi(d,k) = \ko_{X}(d,k)|_{X_0}$. Let $H$ be the pull-back of a general line in $¶^2$ via $\sigma$. Then we have $\ko_{X}(d,k)|_{X_t} \cong \ko_{¶^2}(d)$ for $t \neq 0$. Furthermore $\chi(d,k)|_{¶} \cong \ko_{¶^2}(d-k)$ and $\chi(d,k)|_{\F} \cong \ko_{\F}(dH-(d-k)E)$.

We fix $n-b+1$ general points $p_0,p_1,\ldots,p_{n-b}$ on $¶$ and
$b$ general points $p_{n-b+1},\ldots,p_n$ on $\F$. We define
$\kl_0$ to be the linear sub-system of $\chi(d,k)$ defined by all
divisors of $\chi(d,k)$ having multiplicity at least $m_0$ at
$p_0$ and at least $m$ at the points $p_1,\ldots,p_n$ (write
$\kl_0 = \kl(d,m_0,m^{n-b},m^b)$). We say that $\kl_0$ is obtained
from $\kl=\kl(d,m_0,m^n)$ by an \textit{(k,b)-degeneration}.
$\kl_0$ can be considered as a flat limit on $X_0$ of $\kl$. By
semi-continuity we obtain
    \begin{displaymath}
         \ell_0 = \ell(\kl_0) \geq \ell (\kl).
    \end{displaymath}
  In particular, if $\ell_0 = e(\kl)$ then $\kl$ is non-special.

Now $\kl_0$ restricts on $¶$ to a system $\kl_{¶} =
\kl(d-k,m_0,m^{n-b})$. Furthermore we restrict $\kl_0$ on $\F$ to
$\kl_{\F} = \kl(d,d-k,m^b)$ (the identification we obtain by
blowing down $\kl_{\F}$ to $¶^2$ via $\sigma$). Now we define as
in \cite{CilI} $\kr_{¶}$ to be the linear system on $\kr$ obtained
by restricting $\kl_{¶}$ to $\kr$. We have the following
exact sequence
\begin{displaymath}
          0 \lra \hat{\kl}_{¶} \stackrel{+ L}{\lra} \kl_{¶} \stackrel{|_{L}}{\lra} \kr_{¶} \lra 0.
\end{displaymath}

The kernel system $\hat \kl_{¶}$ consists of all divisors having $L$
as component. So we can identify $\hat \kl_{¶} = \kl(d-k-1,m_0,m^{n-b})$.

We analogously define $\kr_{\F}$ and obtain $\hat \kl_{\F} = \kl(d,d-k+1,m^b)$ (parametrising the divisors in $\kl_{\F}$ which have $E$ as a component).

Let us recall some further abbreviations from \cite{CilI}:
\begin{definitions}

\label{def-dim}

\renewcommand{\labelitemi}{}

   \begin{itemize}

    \item

    \item $v_{¶} = v(\kl_{¶})$,  $v_{\F} = v(\kl_{\F})$,

    \item $\hat{v}_{¶} = v(\hat{\kl}_{¶})$,   $\hat{v}_{\F} = v(\hat{\kl}_{\F})$,

    \item $\ell_{¶}= \ell(\kl_{¶})$,  $\ell_{\F}= \ell(\kl_{\F})$,

   \item $\hat{\ell}_{¶}= \ell(\hat{\kl}_{¶})$, $\hat{\ell}_{\F}= \ell(\hat{\kl}_{\F})$,

  \item  $r_{¶} = \ell_{¶} - \hat{\ell}_{¶} - 1$, the dimension of $\kr_{¶}$,

  \item   $r_{\F} = \ell_{\F} - \hat{\ell}_{\F} - 1$, the dimension of $\kr_{\F}.$

  \end{itemize}

\end{definitions}

In \cite{CilI} it is shown that the associated vector spaces to $\kr_{¶}$ and $\kr_{\F}$ are transversal subspaces of $\Gamma(\kr,\ko_{\kr}(d-k))$.
This leads to the following corollary:

\begin{corollary}[Key-Lemma on $\ell_0$]{We have two cases:}
\label{key-lemma}
\medskip

\begin{enumerate}
 \item If $r_{¶} + r_{\F}$ $\leq$ $d-k-1$, then $\ell_0 = \hat{\ell}_{\mathbb{P}} + \hat{\ell}_{\mathbb{F}} + 1.$

\medskip

 \item If $r_{¶} + r_{\F}$ $\geq$ $d-k-1$, then $\ell_0 = \ell_{\mathbb{P}} + \ell_{\mathbb{F}}-d+k.$

\end{enumerate}

\end{corollary}

A proof can be found in \cite{CilI}

\section{Proof of the Main Theorem}

Before giving the proof let us state two lemmas which are corollaries of the Key-Lemma \ref{key-lemma}. The proof of these is given for an analogous case in \cite{Laf}.

\begin{lemma}[case $v(\kl) \leq -1$]
\label{lemmav<-1}
Let $\kl =\kl(d,m_0,6^n)$ with $v(\kl) \leq -1$. If there are integers $k$
$( k < d)$ and $b$ $( b < n)$ such that a $(k,b)$-degeneration can be
found with the following properties of the restrictions of $\kl_0$
\begin{itemize}
  \item $\kl_{\F}$ and $\kl_{¶}$ are both non-special, and
  \item the kernel systems $\hat\kl_{\F}$ and $\hat \kl_{¶}$ are empty with $\hat v_{¶} \leq v(\kl)$,
\end{itemize}
then $\kl$ is empty.

\end{lemma}

\begin{lemma}[case $v(\kl) \geq -1$]
\label{lemmav>-1}
Let $\kl =\kl(d,m_0,6^n)$ with $v(\kl) \geq -1$. If there are integers $k$ $( k < d)$
and b $( b < n)$ such that a $(k,b)$-degeneration can be found with
\begin{itemize}
  \item $\kl_{\F}$ and $\kl_{¶}$ are both non-special, $v_{¶} \geq -1$, $v_{\F} \geq -1$, and

  \item the kernel systems $\hat \kl_{\F}$ and $\hat \kl_{¶}$ have the property $v(\kl)-1 \geq \hat{\ell}_{¶} + \hat{\ell}_{\F}$,

\end{itemize}

then $\kl$ is non-special.
\end{lemma}

The following three lemmas state parts of the result of the Main Theorem \ref{main-theorem}. We prove them independently later on.

\begin{lemma}[three base points]
 \label{three-points}
  A linear system $\kl(d,m_0,m^n)$ with at most three base points ($n \leq 2$) is special if and only if it is $(-1)$-special.
\end{lemma}

\begin{lemma}[large multiplicities $m_0$ in $p_0$]
 \label{huge-multiplicities}
  Let $d \geq 25$.
  If $m_0 \geq d-9$ then $\kl(d,m_0,6^n)$ is special if and only if
  it is $(-1)$-special.
\end{lemma}

\begin{lemma}[low degrees]
\label{low-degrees}
  If $d \leq 140$ then $\kl(d,m_0,6^n)$ is special if and only if
  it is $(-1)$-special.
\end{lemma}

\begin{proof}[Proof of the Main Theorem \ref{main-theorem}]

  Let $\kl = \kl(d,m_0,6^n)$. By the lemma for large multiplicities
  (\ref{huge-multiplicities}) we can assume that $d \geq m_0 + 10 \geq
  10$.

  Furthermore by the lemma for low degrees (\ref{low-degrees}) the
  statement is true for $d \leq 140$. We can assume $d \geq 141$.
  We continue by induction on $d$ where \ref{low-degrees} can be considered as the base of the induction.

  As all such $\kl$ are not $(-1)$-special we have to show that $\kl$ is non-special. The method is to get the
system $\kl_0$ on the special
  fiber by a degeneration of $\kl$. With Lemmas \ref{lemmav<-1} and \ref{lemmav>-1} we can prove the regularity of $\kl$ if the restrictions of $\kl_0$ to $¶$ and to $\F$ have certain properties. These properties can be achieved as the main conjecture holds for the systems on $¶$ by induction and for the ones on $\F$ by \ref{huge-multiplicities}.

  \medskip

  We perform now a $(5,b)$-degeneration on $\kl$ and get the following systems on the special fiber:
 \medskip

     \begin{tabular}{llll}
       $¶$: &$\kl_{¶}=\kl(d-5,m_0,6^{n-b})$&$\F$:&$\kl_{\F}=\kl(d,d-5,6^b)$\\
         \\
       &$\hat{\kl}_{¶}=\kl(d-6,m_0,6^{n-b})$&&$\hat{\kl}_{\F}=\kl(d,d-4,6^b)$\\
     \end{tabular}
  \medskip

\ul{Step 1 (case $v(\kl) \leq -1$):}

 We want to apply Lemma \ref{lemmav<-1} for the case $v(\kl) \leq -1$.

First of all we need to have $\hat{\kl}_{\F}$ empty. By the lemma
for large multiplicities in $m_0$ (\ref{huge-multiplicities}) we
have that $\hat{\kl}_{\F}$ is non-special if it is
non-$(-1)$-special. Therefore by our classification theorem
\ref{classification} it is sufficient to choose $d < 4b$, i.e., $b
> \frac{d}{4}$. Also we get $\hat v_{\F} \leq -1$, which means this
system is empty.

Next let us find a sufficient condition to get $\hat v_{¶} \leq v(\kl)$. A
computation gives $\hat v_{¶}- v(\kl) = -6d+21b+9$, hence it is sufficient
to have $-6d+21b+9 \leq 0$, that is $b \leq \frac{6d-9}{21}$.

Now we want to find sufficient conditions to have $\kl_{\F}$
non-special. By \ref{huge-multiplicities} this is already the case if we find conditions for $\kl_{\F}$ not to be $(-1)$-special. By Theorem \ref{classification} it is sufficient to force $d > \frac{7b}{2}+3$, that is $b < \frac{2}{7}(d-3)$. As $ \frac{2}{7}(d-3) \leq \frac{6d-9}{21}$,
this new condition on $b$ includes also $\hat v_{¶} \leq v(\kl)$.

In the next step we are searching for a sufficient condition to
get $\kl_{¶}$ non-special. By induction on $d$ $\kl_{¶} =
\kl(d-5,m_0,6^{n-b})$ is special if and only if it is
$(-1)$-special. By our list in Theorem \ref{classification}
we notice that $\kl_{¶}$ is non-$(-1)$-special if we choose
$n-b$ odd as we have assumed that $d-m_0 \geq 10$ and $d \geq
141$.

In the last step we look for a sufficient condition on $b$ to get
$\hat{\kl}_{¶}$ empty. Here we have to be more careful. When
$d-m_0 \geq 11$ we get for the same reasons as in the case of
$\kl_{¶}$ that $\hat{\kl}_{¶}$ is non-special if $n-b$ is
odd. When $d-m_0 =10$ then from Theorem \ref{classification} we
know that $\hat{\ell}_{¶} = -20(n-b) + 5(d-6) - 6$ if $n-b$ is odd.
That means we want this expression to be negative. From
$\hat{\ell}_{¶} \leq -1 \Llra b \leq \frac{1}{4}(7-d)+n$ we get a
sufficient condition on $b$. As by assumption $v(\kl) \leq -1$, we
can conclude that $v(\kl)=11d-21n-45 \leq -1$. Therefore $n \geq
\frac{11d-44}{21}$. That means we can formulate the above
condition on $b$ without $n$ (using a lower bound on $n$) and get
$b \leq \frac{1}{4}(7-d)+\frac {11 d-44}{21} = -
\frac{29}{84}+\frac{23 d}{84}$.

Let us now reformulate all sufficient conditions (separated for the cases $d-m_0 = 10$ and $d-m_0 > 10$) in a compact form:
If $d-m_0 > 10$ we find a $b$ such that we can apply Lemma \ref{lemmav<-1} if
   \begin{displaymath}
      \frac{2}{7} d - \frac{6}{7} - \frac {1}{4}d > 2   \Llra d \geq 81.
   \end{displaymath}
If $d-m_0 = 10$ we find also a $b$ to apply \ref{lemmav<-1} if
   \begin{displaymath}
       - \frac{29}{84}+\frac{23}{84}d - \frac{1}{4} d > 2 \Llra d \geq 99.
   \end{displaymath}

\medskip
\ul{Step 2 (case $v(\kl) \geq -1$):}

We want to use Lemma \ref{lemmav>-1} for the case $v(\kl) \geq -1$. Still all notations are with respect to the above $(5,b)$-degeneration.

In a first step we want to find a sufficient condition on $b$ to get $\kl_{¶}$ non-special. Exactly as in step 1 we get by induction that $\kl_{¶}$ is non-special if we choose $b$ such that $n-b$ is odd, because we assume $d-m_0 \geq 10$ and $d \geq 141$.

Next we want to find sufficient conditions on $b$ to get the system
$\kl_{\F}$ non-special and $v_{\F} \geq -1$. By the lemma for large multiplicities (\ref{huge-multiplicities}) in $m_0$ we have
again as above that $\kl_{\F}$ is non-special if and only if
it is non-$(-1)$-special. We conclude that we get $\kl_{\F}$ non-special if we have $d > \frac{7}{2}b+3$, that is if $b <
\frac{2}{7}d - \frac{6}{7}$, by Theorem \ref{classification}. As
$v_{\F} = 6d-21b-10$ we see that $v_{\F}\geq -1$ which is equivalent to $b \leq
\frac{2}{7}d-\frac{3}{7}$. Therefore the condition for getting
$\kl_{\F}$ non-special gives already that $v_{\F} \geq -1$.

From Theorem \ref{classification} we note again that $b > \frac{d}{4}$ confirms that $\hat{\kl}_{\F}$ is non-special and $\hat{v}_{\F} \leq -1$.

Let us now consider $\hat{\kl}_{¶}$: As above $\hat{\kl}_{¶}$ is
by induction non-special if $n-b$ is odd and $d-m_0 \geq
11$. In the case $d-m_0 \geq 11$ we force also $\hat{v}_{¶} \leq
v(\kl)$, that is $b \leq \frac{6d-9}{21}$. In the case $d-m_0 = 10$ we conclude - exactly as above - that if $n-b$ is odd
$\hat{\kl}_{¶}$ is non-special or $\hat{\ell}_{¶} =
-20(n-b)+5(d-6)-6$. Therefore we force $-20(n-b)+5(d-6)-6 \leq -1$, that means $b \leq \frac{1}{4}(7-d) +n$. As we are in the case $v(\kl) \geq -1$
we have the equation $11d -21n-45 \geq -1$ which means $n \leq
\frac{11d-44}{21}$. It is enough to check the independence of all
conditions on the base points in $\kl$ for the highest possible
number $n$ of points. We fix this $n$ and use a lower bound
$\frac{11d-44}{21}-1$ of it. That means a sufficient condition for
$\hat{\kl}_{¶}$ to be non-special is $b \leq
\frac{1}{4}(7-d)+\frac{11d-44}{21}-1 = \frac{-113+23d}{84}$.

To fulfill all these conditions we need to have $d$ large enough. All
together this gives so far:

If $d-m_0 \geq 11$ we are able to find a sufficient $b$ if
  \begin{displaymath}
     \frac{2}{7}d-\frac{6}{7}- \frac{1}{4}d > 2 \Llra d \geq 81.
  \end{displaymath}
If $d-m_0 = 10$ we are able to find a sufficient $b$ if
  \begin{displaymath}
     \frac{23}{84}d-\frac{113}{84}-\frac{1}{4}d > 2 \Llra d \geq 141.
  \end{displaymath}
In both cases we have that $\kl_{¶}$ and $\kl_{\F}$ are non-special and $v_{\F}\geq -1$. From $\hat{v}_{\F} \leq -1$ and from $v_{¶} = v -1-\hat{v}_{\F}$ we get immediately $v_{¶}\geq -1$. We have $v \geq \hat{v}_{¶}$. As $\hat \kl_{\F}$ and $\hat \kl_{¶}$ are non-special we are able to conclude the following two cases:

If $\hat{v}_{¶} \leq -1$ then
   \begin{displaymath}
     \hat{\ell}_{¶} + \hat{\ell}_{\F} = -2 \leq v - 1,
   \end{displaymath}
and if $\hat{v}_{¶} \geq -1$ then
   \begin{displaymath}
     \hat{\ell}_{¶} + \hat{\ell}_{\F} = \hat{v}_{¶} + \hat{\ell}_{\F} \leq v - 1.
   \end{displaymath}
In both cases we are able to apply Lemma \ref{lemmav>-1} and conclude that $\kl$ is non-special.

\end{proof}

\section{Proof of the Lemmas}

Before starting the proofs we should take some time to explain the use of Quadratic Cremona Transformations for our
purpose. We identify such a transformation with blowing up three general points and
blowing down their connecting lines. Such a transformation is called
to be \textit{based on the three points}. Furthermore one can see by
the blow-up and -down interpretation that a linear system
$\kl(d,m_0,m_1,m_2,m_3,\ldots,m_n)$ is transformed by a Cremona
transformation based on the points $p_0,p_1,p_2$ to a system
$\kl(2d-m_0-m_1-m_2,d-m_1-m_2,d-m_0-m_2,d-m_0-m_1,m_3,\ldots,m_n)$. If all
involved numbers are non-negative (see \cite{CilI}), the dimension and
the virtual dimension of a system $\kl$ do not change under Cremona transformations. In fact a $(-1)$-curve splitting off a system $\kl$ is transformed again into a $(-1)$-curve, which splits off the transformed system. Therefore it is equivalent to examine a system $\kl$ or its Cremona transformed for our purpose. We use suitable sequences of Cremona transformations in the following proofs to obtain systems which are already examined in previous papers.

\begin{proof}[Proof of the lemma of three base points \ref{three-points}]
This can be seen by direct computations with base points $(1:0:0)$, $(0:1:0)$ and $(0:0:1)$. Of course, the statement is also included in the result in \cite{Har}.
\end{proof}

\begin{proof}[Proof of the lemma of large multiplicities $m_0$ in $p_0$ \ref{huge-multiplicities}]
We consider the system $\kl(d,m_0,6^n)$. For the case of $m_0 \geq d-7$ \cite[Proposition 6.2., Corollary 6.3., Proposition 6.4.]{CilI} give a classification of the special systems of this type. Comparing it with our list in \mbox{Theorem \ref{classification}} gives the statement.
Now let $d \geq 25$. The strategy for the proof is to perform a
sequence of Cremona transformations in order to get systems, which can
be examined easier. Furthermore we apply the degeneration
  method again and use again Cremona transformations to prove regularity of some of the obtained systems.

\medskip

\ul{case: $d-m_0 = 8$}

Let $\kl=\kl(d,d-8,6^n)$. We note that if we perform $k$ Cremona transformations, based on $p_0$ and successively on two other base points of multiplicity $6$, we obtain that it is now equivalent to consider the Cremona transformed system (for the strategy see \cite{Laf}):
\begin{displaymath}
   \kl \sim \kl(d-4k,d-8-4k,6^{n-2k},2^{2k})
\end{displaymath}
We set $d-8 = 4t + \epsilon$ with $\epsilon \in \{0,1,2,3\}$. And $n = 2q + \eta$ with $\eta \in \{0,1\}$.

If $t \leq q$ we perform $k=t$ transformations on $\kl(d,d-8,6^n)$ based on $p_0$ and successively two other base points of multiplicity $6$ and obtain
\begin{displaymath}
   \kl \sim \kl(8+\epsilon,\epsilon,6^{n-2t},2^{2t}).
\end{displaymath}
The system on the right hand side is of bounded multiplicity, that means all multiplicities are $\leq 6$. Such systems are special if and only if they are $(-1)$-special by \cite{Yang}.

If $t > q$ we perform $k = q$ transformations on $\kl(d,d-8,6^n)$ again based on $p_0$ and successively two other base points of multiplicity $6$ and obtain
\begin{displaymath}
   \kl \sim \kl(d-4q,d-8-4q,6^{\eta},2^{2q}).
\end{displaymath}
If $\eta = 0$ we are in the case of quasi-homogeneous linear systems of multiplicity $2$, here the main conjecture is true by \cite{CilI}.

If $\eta = 1$ we have to examine systems of the type $\kl = \kl(\delta,\delta-8,6,2^{2q})$ with $\delta = d-4q$. Now let us perform a $(2,b)$-degeneration and get the following systems:

\renewcommand{\labelitemi}{}

\begin{itemize}
 \item $\kl_{¶} = \kl(\delta-2,\delta-8,6,2^{2q-b})$ \space     $\kl_{\F}=\kl(\delta,\delta-2,2^{b})$
 \item $\hat{\kl}_{¶} = \kl(\delta-3,\delta-8,6,2^{2q-b})$ \space    $\hat{\kl}_{\F}=\kl(\delta,\delta-1,2^{b})$

\end{itemize}

If $v(\kl) \leq -1$ we want to apply lemma \ref{lemmav<-1}.

By our classification Theorem \ref{classification} there is no $(-1)$-special system of the type $\kl(d,d-8,6^n)$ if $d \geq 25$. That means we have to show that the system $\kl$ is empty. To use \ref{lemmav<-1} we have again to consider all the systems obtained by the degeneration as in the proof of the main theorem.

In a first step let us consider $\hat{\kl}_{\F}$. As $\hat{\kl}_{\F}$ is a quasi-homogeneous system of multiplicity $m = 2$ we see in \cite{CilI}, that this system is never special. Then $\hat{v}_{\F} = 2 \delta - 3 b$ leads to a sufficient condition to get $\hat{\kl}_{\F}$ empty. This condition is $b \geq \frac{2 \delta +1}{3}$.

In a next step we want to find a sufficient condition to get $\kl_{\F}$ non-special. This is true by \cite{CilI} if $b$ is odd. So let us force $b$ to be odd as a sufficient condition for this case.

Now we consider $\kl_{¶}$. We claim: $\kl_{¶}$ is non-special.

To show the claim we apply at first a Cremona transformation based on the points of multiplicity $\delta-8$, $6$ and on one point of multiplicity $2$. This leads to the following system:
   \begin{displaymath}
   \kl_{¶} \sim \kl(\delta - 4,\delta - 10,4, 2^{2q - b -1}).
   \end{displaymath}
Above we forced $b$ to be odd, therefore we assume $2q-b-1 \geq 2$ (otherwise skip this step) is even. Now we apply successively $\frac{2q-b-1}{2}$ Cremona transformations, based in $p_0$ and two points of multiplicity $2$.
Therefore we see that we have the following equivalence:
  \begin{displaymath}
   \kl_{¶} \sim \kl(\delta-4+2q-b-1,\delta-10+2q-b-1, 4^{2q-b}).
  \end{displaymath}
From $\delta = d-4q \geq 12+\epsilon$ we get by \cite[Theorem 2.1, Theorem 5.2]{Seib} that this system is never special.

Finally we have to consider $\hat{\kl}_{¶}$. Again we claim that $\hat{\kl}_{¶}$ is never special.

 We have by the above assumption that $2q-b$ is odd. At first we split off the line through the points of multiplicity $\delta-8$ and $6$. As the virtual dimension doesn't change we get
  \begin{displaymath}
   \hat{\kl}_{¶} \sim \kl(\delta-4,\delta-9,5,2^{2q-b}).
  \end{displaymath}
Another Cremona transformation based in $p_0$, $p_1$ and one point of multiplicity 2 leads to the equivalence
  \begin{displaymath}
   \hat{\kl}_{¶} \sim \kl(\delta-6,\delta-11,3,2^{2q-b-1}).
  \end{displaymath}
Now as in the case of $\kl_{¶}$ we apply another $\frac{2q-b-1}{2}$ Cremona transformations based in $p_0$ and successively in two points of multiplicity $2$. We end up with the equivalence:
  \begin{displaymath}
   \hat{\kl}_{¶} \sim \kl(\delta-6+\frac{2q-b-1}{2},\delta-11+\frac{2q-b-1}{2}, 3^{2q-b}).
  \end{displaymath}
Now we are able to conclude with \cite{CilI} - as we are in the case of a quasi-homogeneous system of multiplicity 3 - that this system is never special.

To apply \ref{lemmav<-1} we have to find a sufficient condition for $b$ to get $\hat{v}_{¶} \leq -1 $, therefore it is sufficient to have $\hat{v}_{¶} - v(\kl) \leq 0$, which is equivalent to $b \leq \delta$.

All together we find a sufficient $b$ if $\delta - \frac{2 \delta + 1}{3} \geq 2 \Llra \delta \geq 8$. As we have seen above we have already $\delta \geq 12+\epsilon$. This means we can apply Lemma \ref{lemmav<-1} and conclude that $\kl(d,d-8,6^n)$ is empty in the case $v(\kl) \leq -1$.

Now we have to consider the case $v(\kl) \geq -1$. Here we want to apply the Lemma \ref{lemmav>-1}.

As in the case $v(\kl) \leq -1$ we can always find a $b$ such that all
the systems obtained by the above $(2,b)$-degeneration are
non-special. Let us choose such a $b$ like above and then
consider the systems $\kl_{¶}$, $\hat{\kl}_{¶}$, $\kl_{\F}$ and
$\hat{\kl}_{\F}$. From $v_{¶}=v(\kl)-\hat{v}_{\F}-1$, $\hat{v}_{\F} \leq -1$ and $v(\kl) \geq -1$ we conclude $v_{¶} \geq v(\kl) \geq -1$. A direct computation gives $v_{\F} \geq -1$.

As the inequality $\hat{v}_{¶} \leq v(\kl)$ is also fulfilled we get $\hat{\ell}_{¶} \leq v(\kl)$. Therefore we can apply Lemma \ref{lemmav>-1} and conclude that $\kl(d,d-8,6^n)$ is non-special.

\medskip

\ul{case: $d-m_0 = 9$}

Let $\kl=\kl(d,d-9,6^n)$. We note as above that if we perform $k$ Cremona transformations, based on $p_0$ and successively on two other base points of multiplicity $6$, we obtain that:
\begin{displaymath}
   \kl \sim \kl(d-3k,d-9-3k,6^{n-2k},3^{2k})
\end{displaymath}
We set $d-9 = 3t + \epsilon$ with $\epsilon \in \{0,1,2 \}$. And $n = 2q + \eta$ with $\eta \in \{0,1\}$.

If $t \leq q$ we perform $k=t$ transformations on $\kl(d,d-9,6^n)$ based on $m_0$ and successively on two other base points of multiplicity $6$ and obtain
\begin{displaymath}
   \kl \sim \kl(9+\epsilon,\epsilon,6^{n-2t},3^{2t}).
\end{displaymath}
Then the system on the right hand side is of bounded multiplicity,
that means all multiplicities are $\leq 6$. As mentioned above such
systems are special if and only if they are $(-1)$-special by \cite{Yang}.

If $t > q$ we perform $k=q$ transformations on $\kl(d,d-9,6^n)$ and obtain
\begin{displaymath}
   \kl \sim \kl(d-3q,d-9-3q,6^{\eta},3^{2q}).
\end{displaymath}
If $\eta = 0$ we are in the case of quasi-homogeneous linear systems of multiplicity $3$, here the main conjecture is true by \cite{CilI}.

If $\eta = 1$ we have to examine systems of the type $\kl(\delta,\delta-9,6,3^{2q})$ with $\delta = d-3q$. If $\delta < 15$ we are in the case of systems of bounded multiplicity where the main conjecture holds by \cite{Yang}. So we can assume $\delta \geq 15$. Also we can assume $q \geq 1$ (otherwise the statement is clear). Now let us perform a $(3,b)$-degeneration and get the following systems:

\renewcommand{\labelitemi}{}

\begin{itemize}

\item $\kl_{¶} = \kl(\delta-3,\delta-9,6,3^{2q-b})$  \space  $\kl_{\F}=\kl(\delta,\delta-3,3^{b})$
\item $\hat{\kl}_{¶} = \kl(\delta-4,\delta-9,6,3^{2q-b})$ \space   $\hat{\kl}_{\F}=\kl(\delta,\delta-2,3^{b})$

\end{itemize}

\medskip

If $v(\kl) \leq -1$ we again want to apply Lemma \ref{lemmav<-1}. 

\medskip

So let as go through all the systems from the above $(3,b)$-degeneration and search for sufficient conditions on $b$ to apply Lemma \ref{lemmav<-1}.

Let us consider $\hat{\kl}_{\F}$ at first. Here it is sufficient to choose $b > \frac{\delta}{2}$ to get this system non-special by \CilI\space and $\hat{\ell}_{\F} = -1$.

In a next step consider $\kl_{\F}$. By \CilI\space this is non-special if $b$ is odd. 

Then we force (to apply \ref{lemmav<-1}) $\hat{v}_{¶} \leq v$. This is fulfilled if $b \leq \frac{2 \delta -1}{3}$.

Now let us consider $\kl_{¶}$. We claim that this system is never special. To see that let us perform Cremona transformation based on the points of multiplicity $\delta - 9$, 6 and 3. We obtain:
     \begin{displaymath}
   \kl_{¶} \sim \kl(\delta - 6,\delta - 12,3^{2q - b -1}) 
   \end{displaymath}
These systems are always regular by \CilI\space as we have $\delta$ high enough. 

\medskip

A little bit more complicated is the case of $\hat{\kl}_{¶}$. We are searching for a sufficient condition on $b$ to get $\hat{\kl}_{¶}$ empty. We want to show, that $\hat{\kl}_{¶}$ is never special. Then we get the condition simply be choosing $b$ such that $\hat{v}_{¶} \leq -1$ (fulfilled by $\hat{v}_{¶} \leq v$). 

First of all we split off a line through $p_0$, the point of multiplicity $m_0 = \delta -9$, and the point of multiplicity $6$. Therefore we obtain
    \begin{displaymath}
   \hat{\kl}_{¶} \sim \kl(\delta - 5,\delta - 10,5,3^{2q - b}). 
   \end{displaymath}
as the virtual dimension doesn't change in that case. If $2q-b > 0$  applying a further Cremona transformation based in the points of multiplicity $\delta-10$, $5$ and one point of multiplicity $3$ gives 
 \begin{displaymath}
  \hat{\kl}_{¶} \sim \kl(\delta - 8,\delta - 13,2,3^{2q - b -1}). 
  \end{displaymath}
Note that $2q-b-1$ is an even number, as $b$ is odd. We apply now successively Cremona transformations, based on the point $p_0$ and on two other points of multiplicity $3$. It is better again to consider two different cases.

At first assume $\delta - 13 \geq \frac{2q-b-1}{2}$. Then we get
     \begin{displaymath}
  \hat{\kl}_{¶} \sim \kl(\delta - 8- \frac{2q-b-1}{2} ,\delta - 13 - \frac{2q-b-1}{2},2^{2q - b}). 
   \end{displaymath}
By \CilI\space such a system is never special.

Secondly assume $\delta-13 <  \frac{2q-b-1}{2}$. Let $m = \delta - 13$. Then after $m$ such transformations we obtain: 
   \begin{displaymath}
   \hat{\kl}_{¶} \sim \kl(5 ,2^{2m +1},3^{2q-b-1-2m}). 
   \end{displaymath}
Again splitting off a line through two points of multiplicity $3$ (virtual dimension does not change) gives:
 \begin{displaymath}
   \hat{\kl}_{¶} \sim \kl(4,2^{2(m+1) +1},3^{2q-b-1-2(m+1)}). 
  \end{displaymath}
Now two if $2q-b-1-2(m+1) \geq 2$ splitting off lines gives that $\kl_{¶}$ is empty. Secondly if $2q-b-1-2(m+1) = 0$ we have also by \CilI\space  that the system is empty (as $m \geq 2$ by assumption that $\delta \geq 15$).

Taking into account all our conditions on $b$ we require 
 \begin{displaymath}
    \frac{2 \delta -1}{3} - \frac{\delta}{2} > 2 \Llra \delta \geq 15.
 \end{displaymath}
Finally applying Lemma \ref{lemmav<-1} gives that the system $\kl=\kl(\delta,\delta-9,6,3^{2q})$ is empty in the case $v \leq -1$. 
\medskip

Now we have to consider the case $v(\kl) \geq -1$. We want to apply Lemma \ref{lemmav>-1}.

\medskip

As in the case of $v \leq -1$ we get that all the systems $\kl_{\F}$, $\hat{\kl}_{\F}$, $\kl_{¶}$ and $\hat{\kl}_{¶}$ are non-special and $\hat{\ell}_{\F} = -1$ and $\hat{v}_{¶} \leq v$ with a suitable $b$ for the degeneration. 

That means here $\hat{\ell}_{¶} \leq v$ and we can apply Lemma
\ref{lemmav>-1} and conclude that $v(\kl) = \ell$, that means $\kl=\kl(\delta,\delta-9,6,3^{2q})$ is regular.

\medskip

 This finally completes our proof for the case of multiplicities $m_0 = d-8$ and $m_0 = d-9$.

\end{proof}

\begin{proof}[Proof of the lemma of low degrees \ref{low-degrees}]
The main tool for this proof is a computer program which uses
$(5,b)$- and $(6,b)$-degenerations of the plane in order to
prove that certain non-$(-1)$-special systems are
non-special. This algorithm is given by Laface and
Ugaglia in \cite{Laf}. We implemented this algorithm in
\textit{Singular} (see \cite{Sing}). Furthermore to treat the
cases where the degeneration-method fails we implemented a method used by Yang in \cite{Yang}. This method specializes the base points on a line and moves them to infinity. Then it is easier to check if the given conditions on the base points are independent. If this is still the case it proves regularity of a given system.\\
Below we list only the cases in which the program fails. All these
but $10$ cases are solved by ad-hoc methods (mainly Cremona
transformations). The remaining $10$ cases we computed directly
with \textit{Singular} in characteristic $32003$. One can see that
this implies then regularity in characteristic $0$, too.

\fontsize{10}{2ex}

\begin{tabular*}{\textwidth}%
{@{}c@{\extracolsep{\fill}}l@{\extracolsep{\fill}}l@{\extracolsep{\fill}}l@{}}

\\

$d-m_0$ & system & dimension & method \\
           \hline
\\
  $8$ & $\kl=\kl(8,0,6^3)$ & $-1$ & 3-point lemma\\
  $8$ & $\kl=\kl(9,1,6^3)$ & $-1$ & splitting off lines\\
  $14$ & $\kl=\kl(14,0,6^6)$ & $-1$ & Cremona and splitting off lines\\
  $13$ & $\kl=\kl(14,1,6^6)$ & $-1$ & as $\kl(14,0,6^6)$ is empty\\
  $12$ & $\kl=\kl(14,2,6^6)$ & $-1$ & as $\kl(14,0,6^6)$ is empty\\
  $11$ & $\kl=\kl(14,3,6^6)$ & $-1$ & as $\kl(14,0,6^6)$ is empty\\
  $10$ & $\kl=\kl(14,4,6^6)$ & $-1$ & as $\kl(14,0,6^6)$ is empty \\
 $8$ & $\kl=\kl(14,6,6^5)$ & $-1$  & as $\kl(14,0,6^6)$ is empty\\
  $15$ & $\kl=\kl(15,0,6^7)$ & $-1$ & Cremona \\
  $15$ & $\kl=\kl(15,0,6^6)$ & $>-1$ & as $\kl(15,3,6^6)$ is regular \\
  $14$ & $\kl=\kl(15,1,6^6)$ & $>-1$ & as $\kl(15,3,6^6)$ is regular \\
  $13$ & $\kl=\kl(15,2,6^6)$ & $>-1$ & as $\kl(15,3,6^6)$ is regular\\
 $12$ & $\kl=\kl(15,3,6^6)$ & $>-1$ & Cremona and \CilI\\
  $11$ & $\kl=\kl(15,4,6^6)$ & $-1$ & Cremona and splitting off lines\\
  $10$ & $\kl=\kl(15,5,6^6)$ & $-1$ & as $\kl(15,4,6^6)$ is empty\\
  $9$ & $\kl=\kl(15,6,6^6)$ & $-1$ & as $\kl(15,4,6^6)$ is empty\ \\
  $9$ & $\kl=\kl(15,6,6^5)$ & $>-1$ & as $\kl(15,0,6^6)$ is regular\\
  $8$ & $\kl=\kl(15,7,6^5)$ & $>-1$ & Cremona and \CilI\\
  $16$ & $\kl=\kl(16,0,6^8)$ & $-1$ & as $\kl(16,3,6^7)$ is empty\\
  $16$ & $\kl=\kl(16,0,6^7)$ & $>-1$ & as $\kl(16,2,6^7)$ is regular \\
  $15$ & $\kl=\kl(16,1,6^7)$ & $>-1$ &as $\kl(16,2,6^7)$ is regular \\
  $14$ & $\kl=\kl(16,2,6^7)$ & $>-1$ & Cremona and \CilI\\
  $13$ & $\kl=\kl(16,3,6^7)$ & $-1$ & Cremona and splitting off lines\\
  $12$ & $\kl=\kl(16,4,6^7)$ & $-1$ & as $\kl(16,3,6^7)$ is empty\\
  $11$ & $\kl=\kl(16,5,6^7)$ & $-1$ & as $\kl(16,3,6^7)$ is empty \\
  $10$ & $\kl=\kl(16,6,6^7)$ & $-1$ & as $\kl(16,3,6^7)$ is empty\\
  $10$ & $\kl=\kl(16,6,6^6)$ & $>-1$ & as $\kl(16,2,6^7)$ is regular\\
  $9$ & $\kl=\kl(16,7,6^6)$ & $-1$ & Cremona and splitting off lines \\
 $8$ & $\kl=\kl(16,8,6^6)$ & $-1$ & as $\kl(16,7,6^6)$ is empty\\
  $17$ & $\kl=\kl(17,0,6^8)$ & $>-1$ & as $\kl(17,1,6^8)$ is regular\\
  $16$ & $\kl=\kl(17,1,6^8)$ & $>-1$ & Cremona\\
 $15$ & $\kl=\kl(17,2,6^8)$ & $-1$ & Cremona and splitting off lines\\
  $11$ & $\kl=\kl(17,6,6^7)$ & $>-1$ & as $\kl(17,1,6^8)$ is regular\\
  $10$ & $\kl=\kl(17,7,6^7)$ & $-1$ & Cremona and splitting off lines \\
  $9$ & $\kl=\kl(17,8,6^7)$ & $-1$ & as $\kl(17,7,6^7)$ is empty\\
  $8$ & $\kl=\kl(18,10,6^7)$ & $-1$ & Cremona and splitting off lines\\
  $19$ & $\kl=\kl(19,0,6^{10})$ & $-1$ & \CilII \\
  $18$ & $\kl=\kl(19,1,6^{10})$ & $-1$ & as $\kl(19,0,6^{10})$ is empty\\
  $17$ & $\kl=\kl(19,2,6^{10})$ & $-1$ & as $\kl(19,0,6^{10})$ is empty \\
  $15$ & $\kl=\kl(19,4,6^9)$ & $>-1$ & as $\kl(19,5,6^9)$ is regular\\
$14$ & $\kl=\kl(19,5,6^9)$ & $>-1$ & regular by \Yang \\
  $13$ & $\kl=\kl(19,6,6^9)$ & $-1$ & as $\kl(19,0,6^{10})$ is empty\\
  $12$ & $\kl=\kl(19,7,6^9)$ & $-1$ & as $\kl(19,0,6^{10})$ is empty\\
  $9$ & $\kl=\kl(19,10,6^7)$ & $>-1$ & Cremona and \CilII \\
 $8$ & $\kl=\kl(19,11,6^7)$ & $-1$ & Cremona and splitting off lines\\
  $12$ & $\kl=\kl(20,8,6^9)$ & $>-1$ & direct computation with \Sing\space in char$= 32003$\\
\end{tabular*}

\begin{tabular*}{\textwidth}%
{@{}c@{\extracolsep{\fill}}l@{\extracolsep{\fill}}l@{\extracolsep{\fill}}l@{}}

\\

$d-m_0$ & system & dimension & method \\
           \hline
\\ 
  $11$ & $\kl=\kl(20,9,6^9)$ & $-1$ & Cremona and \Laf \\
  $8$ & $\kl=\kl(20,12,6^7)$ & $>-1$ & Cremona and \CilII \\
  $11$ & $\kl=\kl(21,10,6^9)$ & $>-1$ & Cremona and \Yang \\
  $10$ & $\kl=\kl(21,11,6^9)$ & $-1$ & Cremona and \Yang \\
  $9$ & $\kl=\kl(21,12,6^8)$ & $>-1$ & Cremona and \CilII \\
 $8$ & $\kl=\kl(21,13,6^8)$ & $-1$ & Cremona, splitting off lines and \CilII\\
  $22$ & $\kl=\kl(22,0,6^{13})$ & $>-1$ & as $\kl(22,1,6^{13})$ is regular\\
  $21$ & $\kl=\kl(22,1,6^{13})$ & $>-1$ & \Yang \\
  $20$ & $\kl=\kl(22,2,6^{13})$ & $-1$ & \Yang \\
  $19$ & $\kl=\kl(22,3,6^{13})$ & $-1$ & as $\kl(22,2,6^{13})$ is empty\\
  $16$ & $\kl=\kl(22,6,6^{12})$ & $>-1$ & as $\kl(22,1,6^{13})$ is regular \\
  $15$ & $\kl=\kl(22,7,6^{12})$ & $-1$ & direct computation with \Sing\space in char $=32003$\\
  $13$ & $\kl=\kl(22,9,6^{11})$ & $-1$ &  $''$\\
 $11$ & $\kl=\kl(22,11,6^{10})$ & $-1$ & Cremona and \Yang \\
  $10$ & $\kl=\kl(22,12,6^{10})$ & $-1$ & as $\kl(22,11,6^{10})$ is empty \\
  $10$ & $\kl=\kl(22,12,6^{9})$ & $>-1$ & Cremona and \Seib \\
  $9$ & $\kl=\kl(22,13,6^{9})$ & $-1$ & Cremona and splitting off lines\\
  $8$ & $\kl=\kl(22,14,6^{9})$ & $-1$ & as $\kl(22,13,6^9)$ is empty\\
 $12$ & $\kl=\kl(23,11,6^{11})$ & $>-1$ & direct computation with \Sing\space in char $=32003$\\
  $10$ & $\kl=\kl(23,13,6^{10})$ & $-1$ & Cremona and \Seib \\
  $9$ & $\kl=\kl(23,14,6^{9})$ & $>-1$ & Cremona and \CilII \\
 $8$ & $\kl=\kl(23,15,6^{9})$ & $-1$ & Cremona and splitting off lines\\
  $10$ & $\kl=\kl(24,14,6^{10})$ & $>-1$ & Cremona and \CilII \\
  $9$ & $\kl=\kl(24,15,6^{10})$ & $-1$ & Cremona and \Yang \\
  $8$ & $\kl=\kl(24,16,6^{10})$ & $-1$ & as $\kl(24,15,6^{10})$ is empty\\
  $13$ & $\kl=\kl(25,12,6^{13})$ & $-1$ & direct computation with \Sing\space in char $=32003$ \\
  $10$ & $\kl=\kl(25,15,6^{11})$ & $-1$ & Cremona and \Yang \\
  $12$ & $\kl=\kl(26,14,6^{13})$ & $-1$ & direct computation with \Sing\space in char $= 32003$\\
  $10$ & $\kl=\kl(29,19,6^{13})$ & $>-1$ & direct computation with \Sing\space in char $=32003$\\
  $13$ & $\kl=\kl(31,18,6^{17})$ & $-1$ &  $''$\\
  $10$ & $\kl=\kl(31,21,6^{14})$ & $>-1$ & Cremona and \Seib\\
  $10$ & $\kl=\kl(38,28,6^{18})$ & $-1$ & Cremona and \Seib \\
  $13$ & $\kl=\kl(40,27,6^{23})$ & $-1$ & direct computation with \Sing\space in char $=32003$\\
  $10$ & $\kl=\kl(40,30,6^{19})$ & $-1$ & $''$\\
  $10$ & $\kl=\kl(46,36,6^{22})$ & $-1$ & Cremona and \Seib\\

\end{tabular*}

\end{proof}

\fontsize{11}{2ex}

\end{document}